\documentclass[letterpaper, 10 pt, conference]{ieeeconf} 

\IEEEoverridecommandlockouts                             
\usepackage{amsmath} 
\usepackage{amssymb}  
\usepackage{amsfonts}
\usepackage{dsfont}
\usepackage{graphicx}

\makeatletter

\newtheorem{lemma}{Lemma}
\newtheorem{theorem}[lemma]{Theorem}

\newtheorem{example}{Example}
\newtheorem{definition}{Definition}

\newcommand{\bq}{\begin{equation}}
\newcommand{\eq}{\end{equation}}

\newcommand{\X}{\mathcal{X}}

\newcommand{\mR}{\mathbb{R}}

\renewcommand\baselinestretch{0.99}

\newcommand\tc{\tau}

\title{\LARGE \bf
On Energy Conversion in Port-Hamiltonian Systems}

\author{Arjan van der Schaft and Dimitri Jeltsema
\thanks{A.J. van der Schaft is with the Bernoulli Institute for Mathematics, Computer Science and AI, University of Groningen, 
P.O.~Box 407, 9700 AK, The Netherlands. E--mail: a.j.van.der.schaft@rug.nl, D. Jeltsema is with the School of Engineering and Automotive, HAN University of Applied Sciences,
P.O.~Box 2217, 6802 CE Arnhem, The Netherlands. E--mail: d.jeltsema@han.nl.}
}
\begin{document}

\maketitle
\thispagestyle{empty}
\pagestyle{empty}

\begin{abstract}
We study port-Hamiltonian systems with two external ports, and the strategies and limitations for conversion of energy from one port into the other. It turns out that, apart from the cyclo-passivity of port-Hamiltonian systems, this is related to the internal connection structure of port-Hamiltonian systems. A source of motivation for energy conversion is provided by thermodynamics, in particular the Carnot theory of conversion of thermal into mechanical energy. This is extended to general port-Hamiltonian systems which are satisfying structural conditions on their topology; thus generalizing the Carnot-Clausius theory of heat engines. In particular, the operation of Carnot cycles is extended, which is illustrated by the example of a precursor to the Stirling engine, as well as an electromagnetic actuator. Furthermore, alternative energy conversion control schemes such as energy-routers are discussed.
\end{abstract}

\section{INTRODUCTION}
Energy conversion and harvesting are among the most important current engineering problems. This motivates a number of control questions such as the development of a structural theory of (optimal) energy conversion and efficiency, and the development of effective control strategies for achieving this. This paper aims at addressing some of these questions by making use of port-Hamiltonian systems theory (see e.g.  \cite{passivitybook,jeltsema}), which offers a systematic framework for modeling and control of multiphysics systems. This line of research was already initiated in \cite{dim-arjan} by developing the notion of one-port cyclo-passivity, which identifies topological conditions which limit the energy transfer of one port of the system to the other. In the present paper this is extended by defining and analyzing Carnot cycles for one-port cyclo-passive systems, and applying this to a gas-piston system and an electro-mechanical actuator. Furthermore, initial steps are made towards obviating the conditions for one-port cyclo-passivity by establishing direct topological connections by the use of feedback; thus relating to the Duindam-Stramigioli energy router \cite{duindam} as also discussed in \cite{ortega}.  The paper closes by formulating a number of open research problems for the theory of energy conversion for port-Hamiltonian systems.

%%%%%%%%%%%%%%%%%%%%%%%%%%%%%%%%%%%%%%%%%%%%%%%%%%%%%%%%%
%%%%%%%%%%%%%%%%%%%%%%%%%%%%%%%%%%%%%%%%%%%%%%%%%%%%%%%%%
%%%%%%%%%%%%%%%%%%%%%%%%%%%%%%%%%%%%%%%%%%%%%%%%%%%%%%%%%

\section{CYCLO-PASSIVITY}

Before moving on to port-Hamiltonian systems, let us first recall {\it passivity} theory, and especially the less well-known theory of {\it cyclo-passivity}. Consider a standard input-state-output system
\bq
\label{sigma}
\Sigma:\left\{ \,
\begin{aligned}
\dot{x} & =f(x) +G(x)u, \\
y & = h(x), \quad x \in \X, \ u,y \in \mR^m,
\end{aligned}\right.
\eq
on an $n$-dimensional state space manifold $\X$. Cyclo-passivity was coined in \cite{willems1973} generalizing \cite{willems1972}, and was further developed in \cite{hillmoylan1975, hillmoylan1980} with recent extensions in \cite{cyclo}.

\medskip
\hrule
\medskip

\begin{definition}
$\Sigma$ is {\it cyclo-passive} if 
\bq
\label{cyclic}
\int\limits_0^\tc y^\top (t)u(t) dt \geq 0
\eq
for all $\tc\geq 0$ and all $u:[0,\tc] \to \mR^m$ such that $x(\tc)=x(0)$. 
Furthermore, $\Sigma$ is cyclo-passive {\it with respect to} $x^*$ if \eqref{cyclic} holds for all $\tc\geq 0$ and all $u:[0,\tc] \to \mR^m$ such that $x(\tc)=x(0)=x^*$. In case \eqref{cyclic} holds with equality, we speak about cyclo-losslessness (with respect to $x^*$).
\end{definition}

\medskip
\hrule
\medskip

Cyclo-passivity thus means that cyclic motions always require a net nonnegative amount of externally supplied energy. While energy might be released during some sub-interval of the cyclic motion, over the whole time-interval this needs to be compensated by an amount of supplied energy which is at least as large. Clearly, this is a fundamental property of any physical system, where $y^\top u$ is the externally supplied power. The characterization of the set of possible energy {\it storage functions} of a cyclo-passive system is done via the {\it dissipation inequality} \cite{willems1972}.

\medskip
\hrule
\medskip

\begin{definition}
A (possibly extended) function
$S: \X \to -\infty \cup \mR \cup \infty$ satisfies the {\it dissipation inequality} for system $\Sigma$ if
\bq
\label{DIE}
S(x(t_2)) \leq S(x(t_1)) + \int\limits_{t_1}^{t_2} y^\top (t)u(t) dt
\eq
holds for all $t_1 \leq t_2$, all input functions $u:[t_1,t_2] \to \mR^m$, and all initial conditions $x(t_1)$, where $y(t)=h(x(t))$, with $x(t)$ denoting the solution of $\dot{x}=f(x) +g(x)u$ for initial condition $x(t_1)$ and input function $u:[t_1,t_2] \to \mR^m$. If $S$ is an ordinary function $S: \X \to -\infty \cup \mR$ satisfying \eqref{DIE}, then $S$ is called a {\it storage function}.
\end{definition}

\medskip
\hrule
\medskip

By assuming {\it differentiability} of $S$, the dissipation inequality \eqref{DIE} is easily seen \cite{passivitybook} to be equivalent to the {\it differential dissipation inequality} 
\begin{equation*}
\left[\frac{\partial S}{\partial x}(x)\right]^\top \!\! \big(f(x) + G(x)u\big) \leq h^\top(x)u, 
\end{equation*}
or equivalently, 
\begin{equation*}
h(x) = G^\top(x)\frac{\partial S}{\partial x}(x)  \ \text{and} \  \left[\frac{\partial S}{\partial x}(x)\right]^\top \!\!\! f(x) \leq 0. 
\end{equation*}

Obviously, if there exists a storage function for $\Sigma$, then by substituting $x(\tc)=x(0)$ in \eqref{DIE}, it follows that $\Sigma$ is cyclo-passive. A converse is given as follows. Assuming {\it reachability} from a ground-state $x^*$ and {\it controllability} to this same state $x^*$, define the (possibly extended) functions $S_{ac}: \X \to \mR \cup \infty$ and $S_{rc} : \X \to - \infty \cup \mR$ as
\begin{equation}
\begin{aligned}
S_{ac}(x) &:= \!\!\!\! \mathop{\sup_{u,\tc\geq 0}}_{x(0)=x, \, x(\tc)=x^*} - \int\limits_0^\tc y^\top (t)u(t) dt, \\
S_{rc}(x) &:= \!\!\!\!\!\!\! \mathop{\inf_{u,\tc\geq 0}}_{x(-\tc)=x^*,\, x(0)=x}   \ \ \int\limits_{-\tc}^0 y^\top (t)u(t) dt,
\end{aligned}
\end{equation}
which have an obvious interpretation in terms of energy. Indeed, the function $S_{ac}(x)$ denotes the maximal (in fact, supremal) energy that can be recovered from the system at state $x$ while returning the system to its ground-state $x^*$, and the function $S_{rc}(x)$ denotes the minimal (in fact, infimal) energy that is needed in order to transfer the system from the ground-state $x^*$ to $x$. 

It directly follows \cite{willems1972, passivitybook,cyclo} from the properties of `supremum' and `infimum' that $S_{ac}$ and $S_{rc}$ satisfy the dissipation inequality. Hence, if either $S_{ac}$ or $S_{rc}$ are ordinary functions (i.e., not taking values $\pm\infty$), then the system is cyclo-passive. 
The following basic theorem was obtained in \cite{cyclo}, extending the results of \cite{hillmoylan1975, hillmoylan1980}.

\medskip
\hrule
\medskip

\begin{theorem}
\label{cycthm}
Assume $\Sigma$ is reachable from $x^*$ and controllable to $x^*$. Then $\Sigma$ is cyclo-passive with respect to $x^*$ if and only if
\bq
\label{acrc}
S_{ac}(x) \leq S_{rc}(x), \quad \mbox{ for all } x \in \X. 
\eq
Furthermore, if $\Sigma$ is cyclo-passive with respect to $x^*$ then both $S_{ac}$ and $S_{rc}$ are storage functions for $\Sigma$, implying that $\Sigma$ is cyclo-passive. Moreover, if $\Sigma$ is cyclo-passive with respect to $x^*$, then $S_{ac}(x^*) = S_{rc}(x^*)=0$,
and any other storage function $S$ for $\Sigma$ satisfies
\bq
\label{**}
S_{ac}(x) \leq S(x) - S(x^*) \leq S_{rc}(x). 
\eq
If the system is cyclo-lossless from $x^*$ then $S_{ac} = S_{rc}$, and the storage function is {\it unique} up to a constant.
More generally, $\Sigma$ is cyclo-passive with unique (up to a constant) storage function if and only if for every $x \in \X$
\bq
\mathop{\inf_{u,\tc\geq 0}}_{x(0)=x(\tc)=x^*} \, \oint\limits_0^\tc y^\top (t)u(t) dt =0.
\eq
\end{theorem}

\medskip
\hrule
\medskip

Note that the choice of the ground-state $x^*$ in the above theorem is {\it arbitrary}. Indeed, if $\Sigma$ is reachable from and controllable to $x^*$, then so it is from any other ground-state (and, equivalently, for any two states $x_1,x_2$ there exists a trajectory from $x_1$ to $x_2$). Furthermore, by the above theorem, cyclo-passivity with respect to $x^*$ implies cyclo-passivity and thus cyclo-passivity with respect to any other ground state. 
\medskip
\hrule
\medskip

\begin{example}\label{exampleLMI}
A simple example of a passive system which is {\it not} cyclo-lossless, but still has {\it unique} storage function is the mass-spring-damper system
\begin{equation*}
\begin{aligned}
\begin{bmatrix} \dot{q} \\ \dot{p} \end{bmatrix} & =  
\begin{bmatrix} 0 & \frac{1}{m} \\
-k & -\frac{d}{m} \end{bmatrix}
\begin{bmatrix} q \\ p \end{bmatrix}
 + \begin{bmatrix}0 \\1 \end{bmatrix}u, \quad (u = \text{force}) \\[0.5em]
y & =  \begin{bmatrix}0 & \frac{1}{m} \end{bmatrix}
\begin{bmatrix} q \\ p \end{bmatrix} \quad (= \text{velocity)}
\end{aligned}
\end{equation*}
Here, $q$ is the displacement of a linear spring with stiffness $k$, $p$ is the momentum of mass $m$, and ${d > 0}$ is the damping coefficient. Let $x=\text{col}(q,p)$. The dissipation inequality for quadratic storage functions $S(x)=\frac{1}{2}x^\top Qx$ with 
\begin{equation*}
Q=\begin{bmatrix} q_{11} & q_{12} \\ q_{12} & q_{22} \end{bmatrix}
\end{equation*}
 reduces to the linear matrix inequality (LMI)
\begin{equation*}
\begin{array}{c}
\begin{bmatrix}0 & -k \\ \frac{1}{m} & -\frac{d}{m} \end{bmatrix}\!
\begin{bmatrix}q_{11} & q_{12} \\ q_{12} & q_{22} \end{bmatrix} + 
\begin{bmatrix}q_{11} & q_{12} \\ q_{12} & q_{22} \end{bmatrix}\!
\begin{bmatrix} 0 & \frac{1}{m} \\ -k & -\frac{d}{m} \end{bmatrix} \preceq 0,\\[5mm]
\begin{bmatrix}0 & 1 \end{bmatrix}\!
\begin{bmatrix}q_{11} & q_{12} \\ q_{12} & q_{22} \end{bmatrix}=
\begin{bmatrix}0 & \frac{1}{m} \end{bmatrix}.
\end{array}
\end{equation*}
The equality implies $q_{12}=0$ and $q_{22} = \frac{1}{m}$. Substitution in the inequality yields the unique solution $q_{11} = k$, corresponding to a {\it unique} quadratic storage function $S$ which equals the physical energy $H(q,p) = \frac{1}{2m}p^2 + \frac{1}{2}kq^2$. Thus, the mass-spring-damper system is such that for all $x$ we have 
\begin{equation*}
\mathop{\inf_{u,\tc\geq 0}}_{x(0)=x(\tc)=0} \, \oint\limits_0^\tc \frac{p(t)}{m}u(t) dt =0.
\end{equation*}
Indeed, although all possible trajectories involve non-zero energy-dissipation due to the damper $d>0$, leading to 
\begin{equation*}
\oint\limits_0^\tc \frac{p(t)}{m}u(t) dt >0, 
\end{equation*}
this dissipation can be made arbitrarily small.
\end{example}

\medskip
\hrule
\medskip

The stronger notion of {\it passivity} corresponds to the existence of a {\it nonnegative} storage function. 
%In this case no assumptions regarding reachability and controllability need to be made, and 
One considers instead of $S_{ac}$ the (possibly extended) function
\bq
\label{eq:Sa}
S_{a}(x) := \mathop{\sup_{u,\tc\geq 0}}_{x(0)=x}  - \! \int\limits_0^\tc y^\top (t)u(t) dt,
\eq
which is the supremal energy that can be extracted from the system at state $x$; whence the name {\it available storage} for $S_a$. Obviously, $S_a$ is nonnegative, and, again by the property of `supremum', satisfies the dissipation inequality. 

It follows \cite{willems1972,passivitybook} that $\Sigma$ is passive if and only if $S_a$ is an ordinary function, i.e., $S_a(x) < \infty$ for every $x \in \X$. Furthermore, if $\Sigma$ is passive then $S_a$ is a nonnegative storage function satisfying $\inf_x S_a(x)=0$, and all other {\it nonnegative} storage functions $S$ satisfy
\bq
S_a(x) \leq S(x) - \inf_x S(x),\quad x \in \X.
\eq
Obviously, $S_{ac}(x) \leq S_a(x)$ for any $x \in \X$. If $S_a(x^*)=0$, then it directly follows from the dissipation inequality \eqref{DIE} that passivity implies (and in case $\Sigma$ is reachable from $x^*$, is equivalent to)
\bq
\int\limits_0^\tc y^\top (t)u(t) dt \geq 0,
\eq
for all $u:[0,\tc] \to \mR^m, \tc\geq 0$, where $y(t)$ is the corresponding output for initial condition $x^*$. This is sometimes taken as the definition of passivity, -especially in the linear case with $x^*=0$- , and implies that the supplied energy is nonnegative for any motion from the state of minimal energy. However, in a nonlinear context there is often no natural ground state $x^*$ with $S_a(x^*)=0$, and there may even not exist such a ground state, as demonstrated by the following example \cite{cyclo}.

\medskip
\hrule
\medskip

\begin{example}
\label{exa}
Consider the scalar system
\begin{equation*}
\Sigma:\left\{ \,
\begin{aligned}
\dot{x} &= u, \\
y & = e^x.
\end{aligned}\right.
\end{equation*}
Then, the available storage (\ref{eq:Sa}) takes the form
\begin{equation*}
S_a(x) = \mathop{\sup_{u,\tc\geq 0}}_{x(0)=x} \left\{ e^{x(0)} - e^{x(\tc)} \right\} = e^x,
\end{equation*} 
and $\Sigma$ is lossless with nonnegative storage function $e^x$, which is unique up to a constant. Note that $\inf S_a(x)=0$, but there does not exist a finite $x^*$ such that $S_a(x^*)=0$. Obviously, $\Sigma$ is reachable from and controllable to, e.g., $x^*=0$, with resulting $S_{ac}(x)= S_{rc}(x)= e^x - 1$.
\end{example}

\medskip
\hrule
\medskip

%%%%%%%%%%%%%%%%%%%%%%%%%%%%%%%%%%%%%%%%%%%%%%%%%%%%%%%%%
%%%%%%%%%%%%%%%%%%%%%%%%%%%%%%%%%%%%%%%%%%%%%%%%%%%%%%%%%
%%%%%%%%%%%%%%%%%%%%%%%%%%%%%%%%%%%%%%%%%%%%%%%%%%%%%%%%%

\section{ENERGY CONVERSION BY CARNOT-LIKE CYCLES}\label{sec:one-port_pass}

\begin{figure}[h]
\begin{center}
\includegraphics[width=0.6\columnwidth]{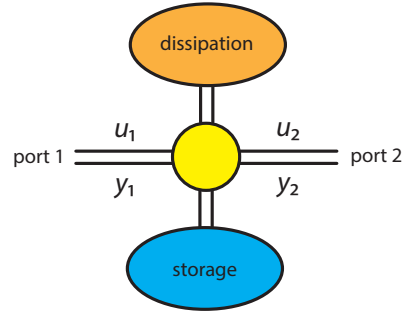}
\caption{Port-Hamiltonian two-port system.}
\label{fig:two-port}
\end{center}
\end{figure}

While the net {\it total} energy supplied to a cyclo-passive system during a cyclic motion is always non-negative, {\it energy conversion} is concerned with the transfer of energy from one port of the system to the other. As we will see, an important aspect in energy conversion is the internal topology of the system. This topology is captured in the port-Hamiltonian formulation of physical systems. A general input-state-output port-Hamiltonian system is given as \cite{passivitybook, jeltsema}
\bq
\label{pH}
\begin{aligned}
\dot{x} & =  J(x)e - \mathcal{R}(x,e) + G(x)u, \quad e= \frac{\partial H}{\partial x}(x), \\[2mm]
y & = G^\top(x)e, \quad x \in \X,
\end{aligned}
\eq
with $n$-dimensional state space $\X$, Hamiltonian $H: \X \to \mR$, skew-symmetric matrix $J(x)=-J^\top(x)$, and dissipation mapping $R$ satisfying $e^\top \mathcal{R} (x,e) \geq 0$ for all $x,e$. By the properties of $J(x)$ and $\mathcal{R}(x,e)$, any port-Hamiltonian system \eqref{pH} satisfies the dissipation inequality
\begin{equation*}
\frac{d}{dt} H(x)= e^\top J(x)e - e^\top \mathcal{R} (x,e) + e^\top G(x)u \leq y^\top u,
\end{equation*}
and thus is {\it cyclo-passive}, with storage function $H$. Conversely, almost any cyclo-passive system can be formulated as a port-Hamiltonian system for some $J$ and $\mathcal{R}$ as above.\footnote{However in systematic systems modeling one would {\it start} from a port-Hamiltonian description with $J$ and $\mathcal{R}$ dictated by the physics of the system.} Thus, by cyclo-passivity theory 
\begin{equation*}
S_{ac}(x) \leq H(x) - H(x^*) \leq S_{rc}(x).
\end{equation*} 
%Also, taking a {\it different} storage function as Hamiltonian will result in a {\it different} $J$ and $R$.
Now, for purposes of energy conversion, consider port-Hamiltonian systems \eqref{pH} with {\it two} ports $(u_1,y_1)$ and $(u_2,y_2)$ as in Fig.~\ref{fig:two-port}. Consequently 
\bq
\frac{d}{dt} H(x) \leq y_1^\top u_1 + y_2^\top u_2.
\eq
The central question addressed in this paper is how we can convert energy which is flowing into the system at port $1$ into energy which is flowing out of the system at port $2$, and to identify the possible obstructions for doing so. Indeed, it turns out that the feasible strategies for energy conversion depend on the structure of the interconnection matrix $J(x)$.

\subsection{Topological Obstructions for Energy Conversion}

Consider the subclass of port-Hamiltonian systems with two ports given as
\bq
\label{form}
\begin{aligned}
\dot{x}_1 & = J_1(x_1,x_2)e_1 -R_1(x_1,x_2)e_1+ G_1u_1,\\[1em]
\dot{x}_2 & = J_2(x_1,x_2)e_2  -R_2(x_1,x_2)e_2 + G_2(x_1,x_2)u_2,\\[0.5em]
y_1 & = G_1^\top e_1, \qquad \quad e_1= \frac{\partial H}{\partial x_1}(x_1,x_2), \\
y_2 & = G_2^\top (x_1,x_2)e_2, \quad  e_2= \frac{\partial H}{\partial x_2}(x_1,x_2),
\end{aligned}
\eq
where $G_1$ is an invertible constant matrix, and the partial Hessian $\frac{\partial^2 H}{\partial x_1^2}(x_1,x_2)$ has full rank everywhere. Because of the first assumption by proper choice of $u_1$ we can ensure that $\dot{x}_1=0$, and thus $x_1=\bar{x}_1$ is constant. It follows that for such $u_1$
\bq
\frac{d}{dt}H(\bar{x}_1,x_2) = \frac{\partial H}{\partial x_2}(x_1,x_2)\dot{x}_2 \leq y_2^\top u_2
\eq
Hence for such input $u_1$ the system is {\it cyclo-passive} at the second port with storage function $H(\bar{x}_1,x_2)$. This directly implies that for such input $u_1$ the system can {\it not} generate energy at port $2$ in a recurrent manner, since by cyclo-passivity the net energy at port $2$ is always flowing {\it into} the system.

Alternatively, by invoking both assumptions, it is clear that, instead of keeping $x_1$ constant, by a proper choice of $u_1$ we can also ensure that $e_1$, or equivalently $y_1$, is constant.

Recall now that the {\it partial Legendre transformation} \cite{lanczos} of $H$ with respect to $x_1$ is given as
\begin{equation*}
\label{legendre}
H^{*}_1(e_1,x_2) = H(x_1,x_2) - e_1^\top x_1, \ e_1= \frac{\partial H}{\partial x_1}(x_1,x_2),
\end{equation*}
where $x_1$ is expressed as a function of $e_1,x_2$ by means of the equation $e_1= \frac{\partial H}{\partial x_1}(x_1,x_2)$ (locally guaranteed by the full rank assumption on the Hessian matrix). 
The following properties of the partial Legendre transformation are well-known \cite{lanczos}
\begin{equation*}
\frac{\partial H^{*}_1}{\partial e_1}(e_1,x_2)= - x_1, \
\frac{\partial H^{*}_1}{\partial x_2}(e_1,x_2)= \frac{\partial H}{\partial x_2}(x_1,x_2). 
\end{equation*}
Hence, with $u_1$ such that $y_1 =\bar{y}_1$ is constant, and, equivalently, $e_1=\bar{e}_1$ is constant, we have that
\begin{equation*}
\begin{aligned}
\frac{d}{dt} H^{*}_1(\bar{e}_1&,x_2) = -x_1^\top \dot{\bar{e}}_1 +  e_2^\top \dot{x}_2 =
e_2^\top J_2(x_1,x_2)e_2\\ 
& -e_2^\top R_2(x_1,x_2)e_2 + e_2^\top G_2(x_1,x_2)u_2 \leq y^\top_2u_2.
\end{aligned}
\end{equation*}
This implies that the system for any constant $y_1$ is {\it cyclo-passive} at the second port $(u_2,y_2)$, with respect to the storage function $H^{*}_1(e_1,x_2)$ with corresponding constant $e_1$. This property was called {\it one-port cyclo-passivity} in \cite{dim-arjan}, and implies that for port-Hamiltonian systems of the form \eqref{form} we cannot convert energy from port $1$ to port $2$ while keeping $y_1$ constant. Furthermore we note that, in view of \eqref{form}, for constant $y_1=\bar{y}_1 =G_1^\top \bar{e}_1$ along all cyclic motions
\begin{equation*}
\begin{aligned}
\int\limits_0^\tc \bar{y}_1^\top u_1(t) dt &= \int\limits_0^\tc \bar{e}_1^\top G_1u_1(t) dt \\[-2.5mm] =\int\limits_0^\tc \bar{e}_1^\top \big[\dot{x}_1&(t) - J_1(x_1,x_2)\bar{e}_1 +R_1(x_1,x_2)\bar{e}_1 \big] dt \\
& \geq\bar{e}_1^\top \left( x_1(\tc) - x_1(0) \right) =0.
\end{aligned}
\end{equation*}
Together this results in the following theorem, extending a main result of \cite{dim-arjan}.

\medskip
\hrule
\medskip

\begin{theorem}
\label{th2}
Consider a port-Hamiltonian system of the form \eqref{form}. Then, for all cyclic motions along which $y_1=\bar{y}_1$ is constant
\bq
\int\limits_0^\tc \bar{y}_1^\top u_1(t) dt \geq 0 \ \  \text{and} \ \int\limits_0^\tc y_2^\top (t) u_2(t) dt \geq 0.
\eq
\end{theorem}

\medskip
\hrule
\medskip

\subsection{Carnot-Cycles and Efficiency}\label{subsec:Carnot}

Within a {\it thermodynamics} context, with $(u_1,y_1)$ corresponding to the {\it thermal} port and $(u_2,y_2)$ corresponding to the {\it mechanical} port, the trajectories for which $x_1$ is constant correspond to {\it isentropic} curves (the entropy $S=x_1$ is constant), also called {\it adiabatics}, while the trajectories for which $e_1=T$ (temperature) are constant are the {\it isothermals}. We will sometimes use this nomenclature for port-Hamiltonian systems as in \eqref{form} as well. 

\begin{figure}[t]
\begin{center}
\includegraphics[width=0.66\columnwidth]{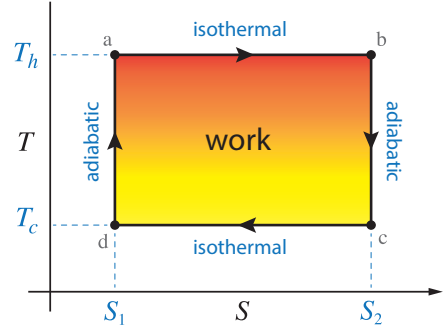}
\caption{Carnot-cycle in the temperature-entropy plane.}
\label{fig:carnot}
\end{center}
\end{figure}

In order to convert energy from port $1$ to $2$ we need by Theorem \ref{th2} more than one value of $e_1$. Within thermodynamics this amounts to the observation by Carnot that at least {\it two} different temperatures $T$ are needed in order to convert thermal energy into mechanical energy (`work') in a recurrent manner. This gives rise to the well-known {\it Carnot cycle}, consisting of an isothermal at high temperature, followed by an adiabatic leading to a lower temperature, another isothermal, and finally an adiabatic which brings the system back to its original state; see Fig.~\ref{fig:carnot}. 

The notion of the Carnot cycle can be extended to general port-Hamiltonian systems of the form \eqref{form} as follows:
\begin{enumerate}
\item On the time-interval $[0,\tau_1]$ consider an `isothermal' with respect to port $1$, corresponding to a constant $e_1=e_1^h$.
Then
\begin{equation*}
\frac{d}{dt} H^*_1(e^h_1,x_2) \leq y_2^\top u_2.
\end{equation*}
\item On the time-interval $[\tau_1, \tau_2]$ consider an `adiabatic' corresponding to a constant $x_1$.
Then
\begin{equation*}
\frac{d}{dt} H(x_1,x_2) \leq y_2^\top u_2.
\end{equation*}
\item On the time-interval $[\tau_2, \tau_3]$ consider an `isothermal'  corresponding to a constant $e_1=e_1^c$.
Then
\begin{equation*}
\frac{d}{dt} H^*_1(e^c_1,x_2) \leq y_2^\top u_2.
\end{equation*}
\item Finally, on the time-interval $[\tau_3, \tc]$ consider an `adiabatic' corresponding to a constant $x_1$. Then
\begin{equation*}
\frac{d}{dt} H(x_1,x_2) \leq y_2^\top u_2.
\end{equation*}
\end{enumerate}
Note that on both isothermals: 
\begin{equation*}
\begin{aligned}
H(x_1,x_2) &= H^*_1(e^h_1,x_2) + e^h_1x_1, \\
H(x_1,x_2) &= H^*_1(e^c_1,x_2) + e^c_1x_1,
\end{aligned}
\end{equation*}
respectively. Since the total process is a cycle, i.e., $x(\tc)=x(0)$, addition of these inequalities yields
\begin{equation*}
\begin{aligned}
0 &= H(x(\tc))  -H(x(0)) \\ 
& \qquad \leq \int\limits_0^\tc y_2^\top (t) u_2(t) dt + e_1^h\Delta^h x_1 + e_1^c\Delta^c x_1,
\end{aligned}
\end{equation*}
where $\Delta^h x_1$ and $\Delta^c x_1$ are the changes in $x_1$ during the isothermals on $[0,\tau_1]$ and $[\tau_3, \tc]$, respectively. Thus, the total energy delivered to the environment via port $2$ satisfies
\bq\label{eq:total_energy}
 -\int\limits_0^\tc y_2^\top (t) u_2(t) dt \leq e_1^h\Delta^h x_1 + e_1^c\Delta^c x_1,
\eq
with equality in case there is no dissipation, i.e., when $R_1$ and $R_2$ are zero. Note that the right-hand side of (\ref{eq:total_energy}) is the total energy provided at port $1$ to the system during the two isothermals. 
Because of $x(\tc)=x(0)$ necessarily $\Delta^h x_1 + \Delta^c x_1=0$. Hence the right-hand side of (\ref{eq:total_energy}) can be also written as $\left(e_1^h - e_1^c \right)\Delta^h x_1$. Assuming that both $e_1^h$ and $e_1^c$ are positive, then this expression is positive if $e_1^h > e_1^c$ and $\Delta^h x_1$. (This is the case considered in thermodynamics where $e_1^h$ and $e_1^c$ are the temperatures of the hot and cold reservoirs.) 
Furthermore, the {\it efficiency} of the Carnot cycle (in case $R_1$ and $R_2$ are zero) is given as
\begin{equation}\label{eq:Carnot-eff}
\frac{e_1^h\Delta^h x_1 + e_1^c\Delta^c x_1}{e_1^h\Delta^h x_1} = 1 - \frac{e_1^c}{e_1^h},
\end{equation}
recovering the well-known expression from reversible thermodynamics \cite{fermi, kondepudi}. 

\medskip
\hrule
\medskip

\begin{figure}[b]
\begin{center}
\includegraphics[width=0.6\columnwidth]{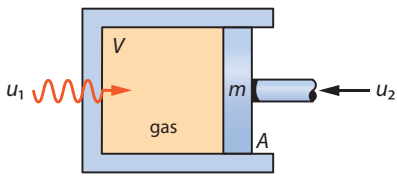}
\caption{Gas-piston system with heat port.}
\label{fig:piston}
\end{center}
\end{figure}

\begin{example}
Consider a gas-piston system with incoming entropy flow (due to incoming heat) as in Fig.~\ref{fig:piston}. The system can be considered as a precursor to the Stirling engine, cf. \cite{stirling}. Its port-Hamiltonian description is given by
\begin{equation}
\label{gaspiston}
\begin{bmatrix}
\dot{S}\\
\dot{V}\\
\dot{\pi}
\end{bmatrix} = 
\begin{bmatrix}
0 & 0  & 0\\
0 & 0  & A\\
0 & -A & 0
\end{bmatrix}\!\!
\begin{bmatrix}
\dfrac{\partial H}{\partial S}\\[0.76em]
\dfrac{\partial H}{\partial V}\\[0.76em]
\dfrac{\partial H}{\partial \pi}\\
\end{bmatrix} + 
\begin{bmatrix}
1 & 0\\
0 & 0\\
0 & 1
\end{bmatrix}\!
\begin{bmatrix}
u_1\\
u_2
\end{bmatrix},
\end{equation}
with Hamiltonian $H(S,V,\pi) = U(S,V) + \frac{1}{2m}\pi^2$ and corresponding natural outputs
\begin{align*}
y_1 &= \frac{\partial U}{\partial S} \ (=T),\\
y_2 &= \frac{\pi}{m}.
\end{align*}
Here, $S$ represents entropy, $T$ temperature, $V$ the volume of the gas enclosed by the piston, $A$ the area of the piston, and $\pi$ the momentum of the piston with mass $m$. Finally, $U(S,V)$ is the internal energy of the gas, $u_1$ the entropy flow, $y_1$ the temperature (so that $y_1u_1$ is the heat flow entering the gas), $u_2$ is the external force on the piston, and $y_2$ its velocity (so that $-y_2u_2$ is the mechanical work done {\it by} the system).
This is a port-Hamiltonian system of the form \eqref{form} with $R_1=R_2=0$. Thus for $u_1$ such that either $S$ is constant (isentropic curves), or $y_1=T$ is constant (isothermals), the system is cyclo-passive at its mechanical port. In order to perform mechanical work in a cyclic manner one needs at least two temperatures $T_h$ (hot) and $T_c$ (cold). 

By using a Carnot cycle consisting of an isothermal at $T_h$, adiabatic, isothermal at $T_c$, and finally an adiabatic to return to the original state $(S,V,\pi)$ the mechanical work delivered to the surrounding is
\bq
\label{stirling}
- \oint y_2(t)u_2(t) dt = \oint y_1(t)u_1(t) dt = Q_h + Q_c,
\eq
where $Q_h$ is the heat entering the system during the isothermal at temperature $T_h$, and $-Q_c$ is the heat leaving the system during the isothermal at temperature $T_c$. 
\end{example}

\medskip
\hrule
\medskip

\begin{example}
Consider an electromagnetic actuator as depicted in Fig.~\ref{fig:mag_actuator}. Such device converts electrical energy into mechanical motion (work). 

\begin{figure}[h]
\begin{center}
\includegraphics[width=0.75\columnwidth]{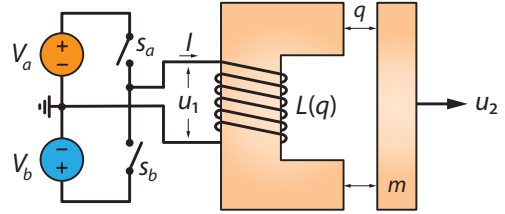}
\caption{Electromagnetic actuator. [Note that in practice a diode must be placed across the inductor terminals.]
}
\label{fig:mag_actuator}
\end{center}
\end{figure}

The port-Hamiltonian description (with $q$ and $p$ the displacement and the momentum of the armature $m$, respectively, and $\varphi$ the magnetic flux of the coil) is given as
  \begin{equation}
  \label{eq:mag_actuator}
    \begin{aligned}
      \begin{bmatrix}\dot \varphi \\ \dot q \\ \dot p 
      \end{bmatrix}  & = 
      \begin{bmatrix}0 & 0 & 0 \\  0 & 0 & 1 \\  0 & -1 &
        0 \end{bmatrix}\!\!
      \begin{bmatrix} 
        \dfrac{\partial H}{\partial \varphi} \\[0.76em] 
        \dfrac{\partial H}{\partial q} \\[0.76em] 
        \dfrac{\partial
          H}{\partial p} 
      \end{bmatrix} + 
      \begin{bmatrix} 1 & 0 \\ 0 & 0\\ 0 &1 \end{bmatrix}\!
      \begin{bmatrix}u_1 \\ u_2 \end{bmatrix}, \\ 
      y_1 & = \dfrac{\partial H}{\partial \varphi} \ (=I), \\
        y_2 & = \dfrac{\partial H}{\partial p},  
    \end{aligned}
  \end{equation}
  where $u_1$ is the supplied voltage and $y_1$ the associated current, while $u_2$ is a mechanical force and $y_2$ the velocity of the armature. This is a port-Hamiltonian system of the form \eqref{form}, with the mechanical and electromagnetic part of the system coupled via the Hamiltonian
\begin{equation*}
H(\varphi,q,p) = \frac{\varphi^2}{2L(q)} + \frac{p^2}{2m},
\end{equation*}
where the electrical inductance ${L(q)>0}$ depends inversely on the horizontal mechanical displacement. 
    
Note that the system (\ref{eq:mag_actuator}) has the same structure as the gas-piston system \eqref{gaspiston}. Thus, (\ref{eq:mag_actuator}) is one-port cyclo-passive at the mechanical port, and a similar Carnot cycle can be used in order to convert energy from the electrical port to the mechanical port. Let $V_a > V_b$. The four processes described in Subsection \ref{subsec:Carnot} are represented as follows:
\begin{enumerate}
\item Switch $s_a$ is closed and $s_b$ is open. The coil absorbs an amount of energy, say $E_a$, from the upper source $V_a$, while the current ${I = I_a}$ remains constant. This results in a decreasing $q$ and a certain amount of external work done by the armature (and port 2). 
\item The coil is disconnected from the voltage source (both switches open) and undergoes an `adiabatic' transformation in which the magnetic flux remains constant and the current $I$ decreases from $I_a$ to $I_b$ as the armature recedes further to the coil. 
\item Switch $s_b$ is now closed and ${I = I_b}$. The armature is turning back to its initial position with energy, say $E_b$, being transferred to the lower voltage source $V_b$. 
\item Both switches open and ${I=I_b}$. The armature displaces to change the inductance $L(q)$ until the current equals ${I = I_a}$. The magnetic flux remains constant.
\end{enumerate}
Now, since the energy 
\begin{equation*}
E_a = \int\limits_{\varphi_a}^{\varphi_b} I_a d \varphi = I_a (\varphi_b - \varphi_a)
\end{equation*}
and, likewise, $E_b = -I_b (\varphi_b - \varphi_a)$, we get that
\begin{equation*}
\frac{E_b}{E_a} = -\frac{I_b}{I_a}.
\end{equation*}
Hence, the efficiency of the system in Fig.~\ref{fig:mag_actuator} equals
\begin{equation*}
\frac{E_a+E_b}{E_a} = 1 + \frac{E_b}{E_a} = 1 - \frac{I_b}{I_a},
\end{equation*}
which is analogous to the efficiency of the Carnot-cycle (\ref{eq:Carnot-eff}). 
\end{example}  
  
\medskip
\hrule
\medskip
  
%%%%%%%%%%%%%%%%%%%%%%%%%%%%%%%%%%%%%%%%%%%%%%%%%%%%%%%%%
%%%%%%%%%%%%%%%%%%%%%%%%%%%%%%%%%%%%%%%%%%%%%%%%%%%%%%%%%
%%%%%%%%%%%%%%%%%%%%%%%%%%%%%%%%%%%%%%%%%%%%%%%%%%%%%%%%%

\section{ENERGY CONVERSION BY DIRECT INTERCONNECTION AND ENERGY-ROUTERS}

The essential assumption in \eqref{form} is the absence of off-diagonal blocks in the $J$-matrix. Since $J$ captures the `topology' of the system, this means that in port-Hamiltonian systems of the form \eqref{form} there is no direct topological connection between the dynamics of the sub-vectors $x_1,x_2$. Instead, the interaction between the two subsystems is via the Hamiltonian $H(x_1,x_2)$ (and to a lesser extent through possible dependence of $J_1$ and $J_2$ on $x_1,x_2$, and similarly for $\mathcal{R}_1, \mathcal{R}_2$). In this section, we consider an example where we {\it do} have non-zero off-diagonal terms (the heat exchanger), and two examples where off-diagonal elements are {\it introduced} through the use of feedback.

\subsection{Heat Exchanger}

\begin{figure}[h]
\begin{center}
\includegraphics[width=0.86\columnwidth]{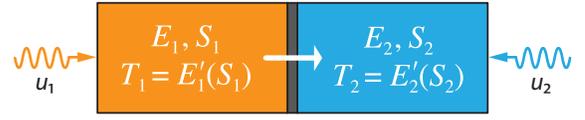}
\caption{Heat exchanger with temperatures $T_1$ and $T_2$.}
\label{fig:hex}
\end{center}
\end{figure}

Consider a heat exchanger comprising two heat reservoirs connected via a conducting wall; see Fig.~\ref{fig:hex}.

\noindent The dynamics is given as
\begin{equation}\label{heatexchanger}
\begin{aligned}
\begin{bmatrix} \dot{S}_1 \\ \dot{S}_2 \end{bmatrix} & =
\begin{bmatrix} 
0 & \!\!\!\!\lambda \frac{E_2' - E_1'}{E_1'E_2'} \\
\lambda \frac{E_1' - E_2'}{E_1'E_2'} & 0 
\end{bmatrix}\!\!
\begin{bmatrix} \dfrac{\partial E}{\partial S_1} \\[0.66em] \dfrac{\partial E}{\partial S_2} \end{bmatrix} 
+  \begin{bmatrix} 1 & 0 \\ 0 & 1 \end{bmatrix} \!\begin{bmatrix} u_1 \\ u_2 \end{bmatrix},\\
y_1 & = E_1'(S_1),\quad y_2 = E_2'(S_2),
\end{aligned}
\end{equation}
where $\lambda$ is Fourier's conduction coefficient, $E_1(S_1)$ and $E_2(S_2)$ are the (thermal) energies of heat reservoirs $1$ and $2$ with temperatures $E_1'(S_1)$ and $E_2'(S_2)$, and 
\begin{equation*}
E(S_1,S_2)=E_1(S_1) + E_2(S_2)
\end{equation*}
represents the total energy of the heat exchanger. Furthermore, $u_1$ and $u_2$ are the external entropy flows entering the reservoirs $1$ and $2$, so that $y_1u_1$ and $y_2u_2$ are the external heat flows entering the two reservoirs. The equations \eqref{heatexchanger} constitute a {\it quasi} port-Hamiltonian system.\footnote{`Quasi' since the $J$-matrix does not directly depend on the state variables $(S_1,S_2)$, but {\it through} the temperatures $T_1=E_1'(S_1)$ and  $T_2=E_2'(S_2)$.}

Although \eqref{heatexchanger} is a thermodynamic system, it is {\it not} in the form \eqref{form}, and, indeed, heat flow occurs from reservoir $1$ to reservoir $2$ for constant $T_1$ as long as $T_1 > T_2$. As an additional property we note that the total entropy $S_1+S_2$ is always non-decreasing.

\subsection{Energy Transfer by Energy-Routing}

Consider two cyclo-lossless port-Hamiltonian systems
\begin{equation*}
\Sigma_i:\left\{
\begin{aligned}
\dot{x}_i &= J_i(x_i) \dfrac{\partial H_i}{\partial x_i}(x_i) + g_i(x_i)u_i, \\
y_i &=  g_i^\top (x_i)\dfrac{\partial H_i}{\partial x_i}(x_i),
\end{aligned}\right.
\end{equation*}
with $i\in\{1,2\}$. The two systems may be coupled to each other by the output feedback (an elementary version of the Duindam-Stramigioli energy-router \cite{duindam, ortega})
\begin{equation*}
\begin{bmatrix} u_1 \\u_2 \end{bmatrix}
=
\begin{bmatrix} 0 & -y_1y_2^\top \\ y_2y_1^\top & 0 \end{bmatrix}\!
\begin{bmatrix} y_1 \\y_2 \end{bmatrix} + \begin{bmatrix} v_1 \\v_2 \end{bmatrix},
\end{equation*}
with new inputs $v_1$ and $v_2$. Then, the closed-loop system $\Sigma_1 \circ \Sigma_2$ is a lossless port-Hamiltonian system with Hamiltonian
\begin{equation*}
H(x_1,x_2)=H_1(x_1) + H_2(x_2),
\end{equation*}
but clearly {\it not} of the form \eqref{form}. Due to the special form of the output feedback, we have that
\begin{equation*}
\begin{aligned}
\frac{d}{dt}H_1 &= - ||y_1||^2 ||y_2||^2 + y_1^\top v_1, \\
\frac{d}{dt}H_2 &=  ||y_1||^2 ||y_2||^2 + y_2^\top v_2.
\end{aligned}
\end{equation*}
Hence, if energy is pumped into $\Sigma_1 \circ \Sigma_2$ through port $1$, then the stored energy in $\Sigma_2$ will increase---enabling its release through port $2$.

\subsection{Energy Transfer Using IDA-PBC}

Consider again the electromagnetic actuator (\ref{eq:mag_actuator}). In order to establish an effective topological coupling between the electrical and mechanical subsystems, we can try to enforce a coupling between the magnetic flux and the momentum; thus realizing a desired interconnection matrix of the form
\begin{equation}\label{eq:Jd_plunger}
J_d = \begin{bmatrix}
0 & 0 & \alpha\\  
0 & 0 & 1 \\  
-\alpha & -1 & 0 
\end{bmatrix},
\end{equation}
where $\alpha$ is a constant or function to be defined. 

The easiest way to accomplish this is when the second input is available for direct manipulation. Indeed, introducing the feedback 
\begin{equation*}
\begin{bmatrix}
u_1\\
u_2
\end{bmatrix} = 
\begin{bmatrix}
0 & \alpha\\
-\alpha & 0
\end{bmatrix}\!
\begin{bmatrix}
y_1\\
y_2
\end{bmatrix}
+\begin{bmatrix}
v_1\\
v_2
\end{bmatrix},
\end{equation*}
with new inputs $v_1$ and $v_2$, readily replaces the $J$-matrix in (\ref{eq:mag_actuator}) by the interconnection matrix (\ref{eq:Jd_plunger}). However, usually only the first input is available for control. 

On the other hand, invoking the IDA-PBC methodology \cite{CSM2001}, and selecting a state feedback $u_1 = \alpha(\varphi,q)\frac{p}{m} + v_1$, with
\begin{equation*}
\alpha(\varphi,q) = \frac{1}{4}\frac{L'(q)}{L(q)}\varphi,
\end{equation*}
and $u_2=v_2$, yields (see the Appendix for details)
\begin{equation*}
    \begin{aligned}
     \begin{bmatrix}
     \dot\varphi\\
     \dot q\\
     \dot p
     \end{bmatrix} =
      \begin{bmatrix}
      0 & 0 & \alpha(\varphi,q) \\  
      0 & 0 & 1 \\  
      -\alpha(\varphi,q) & -1 & 0 \end{bmatrix}\!\!
      \begin{bmatrix} 
        \dfrac{\partial H_d}{\partial \varphi} \\[0.76em] 
        \dfrac{\partial H_d}{\partial q} \\[0.76em] 
        \dfrac{\partial H_d}{\partial p} 
      \end{bmatrix}  + 
      \begin{bmatrix}
      1 & 0\\
      0 & 0\\
      0 & 1
      \end{bmatrix}\!
      \begin{bmatrix}
      v_1\\
      v_2
      \end{bmatrix},
    \end{aligned}
  \end{equation*}
in which the total stored energy is modified into
\begin{equation*}
H_d(\varphi,q,p)=\frac{\varphi^2}{L(q)} + \frac{p^2}{2m}.
\end{equation*}
This (partially) controlled system is clearly {\it not} in the form \eqref{form} anymore, and a direct interaction between the mechanical and magnetic subsystem is achieved through the lossless modulation $\alpha(\varphi,q)$. Moreover, the magnetic energy storage is doubled, enabling an increased energy conversion rate. 

%%%%%%%%%%%%%%%%%%%%%%%%%%%%%%%%%%%%%%%%%%%%%%%%%%%%%%%%%
%%%%%%%%%%%%%%%%%%%%%%%%%%%%%%%%%%%%%%%%%%%%%%%%%%%%%%%%%
%%%%%%%%%%%%%%%%%%%%%%%%%%%%%%%%%%%%%%%%%%%%%%%%%%%%%%%%%

\section{CONCLUSIONS AND OUTLOOK}

Energy conversion in physical systems is related to the structure of the $J$-matrix in its port-Hamiltonian formulation. In particular, if the off-diagonal blocks of $J$ are zero we generally have one-port cyclo-passivity \cite{dim-arjan}, and need more than one value of the output at port $1$ in order to convert energy from port $1$ to port $2$. This may be achieved by a direct generalization of the Carnot cycle for heat engines. Another strategy may be to use feedback in order to introduce non-zero off-diagonal blocks. 

Among the many open problems we mention:
\begin{itemize}
\item Consider a port-Hamiltonian system in the form \eqref{form}. Can we generalize Carnot cycles to more general cyclic processes, while retaining a notion of efficiency?
\item How to formulate {\it optimal} energy conversion from one port to another?
\item The physical energy $H$ of any port-Hamiltonian system is somewhere `in between' $S_{ac}$ and $S_{rc}$. If we choose another storage function this will generally lead to a {\it different} $J$-matrix (as well as different dissipation structure $\mathcal{R}$). How to exploit this for energy conversion?
\end{itemize}

%%%%%%%%%%%%%%%%%%%%%%%%%%%%%%%%%%%%%%%%%%%%%%%%%%%%%%%%%
%%%%%%%%%%%%%%%%%%%%%%%%%%%%%%%%%%%%%%%%%%%%%%%%%%%%%%%%%
%%%%%%%%%%%%%%%%%%%%%%%%%%%%%%%%%%%%%%%%%%%%%%%%%%%%%%%%%

%%%%%%%%%%%%%%%%%%%%%%%%%%%%%%%%%%%%%%%%%%%%%%%%%%%%%%%%%
%%%%%%%%%%%%%%%%%%%%%%%%%%%%%%%%%%%%%%%%%%%%%%%%%%%%%%%%%
%%%%%%%%%%%%%%%%%%%%%%%%%%%%%%%%%%%%%%%%%%%%%%%%%%%%%%%%%

\section*{APPENDIX}

\small

The IDA-PBC methodology \cite{CSM2001} entails to solving the so-called {\it matching} equation
\begin{equation*}
J_d(x)\frac{\partial H_a}{\partial x} = -J_a(x)\frac{\partial H}{\partial x} + G(x)u,
\end{equation*}
with $J_a(x) = J_d(x) - J(x)$. 
For the system (\ref{eq:mag_actuator}), desired interconnection structure (\ref{eq:Jd_plunger}), and $u_1 = \beta(\cdot) + v_1$, this amounts to solving
\begin{equation*}
\begin{aligned} 
\alpha(\cdot)\frac{\partial H_a}{\partial p} &= - \alpha(\cdot) \frac{\partial H}{\partial p} + \beta(\cdot),\\     
\frac{\partial H_a}{\partial p} &=0,\\
-\alpha(\cdot) \frac{\partial H_a}{\partial \varphi} - \frac{\partial H_a}{\partial q} &= \alpha(\cdot) \frac{\partial H}{\partial \varphi}.
\end{aligned}
\end{equation*}
Hence, selecting $H_a(\varphi,q) = \frac{1}{2L(q)}\varphi^2$, we have that 
\begin{equation*}
\frac{1}{2}\frac{L'(q)}{L^2(q)}\varphi^2 = 2 \alpha(\cdot) \frac{\varphi}{L(q)} \Rightarrow \alpha(\varphi,q) = \frac{1}{4}\frac{L'(q)}{L(q)}\varphi,
\end{equation*}
together with the full state feedback control
\begin{equation*}
\beta(\varphi,q,p) = \alpha(\varphi,q) \frac{\partial H}{\partial p},
\end{equation*}
and energy storage $H_d(\varphi,q,p) = H(\varphi,q,p) + H_a(\varphi,q)$. 
\end{document}